\documentclass[12pt]{amsart}
\usepackage{
  fullpage,
  graphicx,amssymb,epic,eepic,cancel,
  picins,
color,hyphenat,
  multicol, 
slashbox, 
}
\usepackage[all]{xy}

\newtheorem{theorem}{Theorem}[section]
\newtheorem{lemma}[theorem]{Lemma}

\begin{document}
\newdimen\captionwidth\captionwidth=\hsize

\title{Pentagon and Hexagon Equations Following Furusho}

\author{Dror Bar-Natan and Zsuzsanna Dancso}
\address{
  Department of Mathematics\\
  University of Toronto\\
  Toronto Ontario M5S 2E4\\
  Canada
}
\email{drorbn@math.toronto.edu, zsuzsi@math.toronto.edu}
\urladdr{http://www.math.toronto.edu/drorbn, http://www.math.toronto.edu/zsuzsi}

\date{September 25th, 2010}


\begin{abstract}
In \cite{F} H. Furusho proves the beautiful result that of the
three defining equations for associators, the pentagon implies
the two hexagons (see also \cite{W}). In this note we present a simpler proof for this
theorem (although our paper is less dense, and hence only slightly shorter). 
In particular, we package the use of algebraic geometry and
Groethendieck-Teichmuller groups into a useful and previously known
principle, and, less significantly, we eliminate the use of
spherical braids.
\end{abstract}

\maketitle
\tableofcontents

\section{Introduction}\label{int}

Associators are useful and intricate gadgets that were first
introduced and studied by Drinfel'd in \cite{Dr0} and \cite{Dr1}.
The theory was later put in the context of parenthesized 
(a.k.a. non-associative) braids by \cite{LM} and \cite{BN}. 
Associators arise in several different areas of mathematics,
and thus constructing an associator is a task faced by many. 
Unfortunately, it
is a very difficult task: no closed formulas are known
at present.

Associators are essentially the solutions to three equations,
called the pentagon and the positive and negative hexagons, which live in complicated
diagrammatic spaces. Finding an associator amounts to 
finding a solution to this system of equations. Furusho's
result states that the last two of these three equations are
superfluous: a solution to the pentagon will automatically
be a solution to both hexagons. This statement is 
quite surprising, and thus we felt that a simpler proof would be
of value.

The paper is organized as follows: we first review some definitions, 
then present the main tool or ``extension principle'' for the simplified proof: a theorem from 
\cite{Dr1} and \cite{BN}, and also a simpler standard fact
which we call the ``linearization principle''. We then 
prove the theorem modulo a ``Lie-algebraic'' version which we call the ``Main 
Lemma'', followed by the proof of the Main Lemma, and finally a side
note on how one of the algebraic maps used in the proof arises from 
topology. 

The part of the proof which depends on the extension principle is 
significantly different from that of \cite{F} (and \cite{W}),
and significantly simpler: Furusho's proof uses algebraic geometry and 
Groethendieck-Teichmuller groups to go from the Lie algebra version to
the group statement. 

The proof of the main lemma is essentially Furusho's
proof, except for the elimination of the use of spherical braids. 
Using
spherical braids makes the statement of the Main Lemma look very symmetric, 
so in some sense prettier. However, in our opinion, it makes the combinatorial 
argument of the proof less transparent: it is the break of symmetry
which helps the reader figure out how someone might have discovered the proof.
Also, we save the trouble of having to pass back and forth between 
the spherical and regular chord-diagram spaces. The elimination of spherical braids
is easy: the space of chord diagrams of pure spherical 
$5$-braids and that of regular pure $4$-braids only differ in taking a quotient by
the center of the latter.

\subsection{Acknowledgements}
The authors would like to thank the referee for the thorough reading
of the paper and detailed helpful suggestions, and Stavros Garoufalidis
for an enlightening discussion of the proof of the Main Lemma.

\section{Definitions}
Algebraically, the space of chord diagrams of pure $n$-braids,
$\mathcal{A}_n$, is defined to be the graded completion of
the following non-associative algebra, where generation
is understood over a fixed field $k$:

\[ \left\langle t_{ij}: 1 \leq i\neq j \leq n
  \left|\ \parbox{3.25in}{
    $t_{ij}=t_{ji}$, \newline
    $[t_{ij},t_{kl}]=0$ when $i,j,k,l$ are distinct (L), \newline
    $[t_{ij}+t_{ik},t_{jk}]=0$ when $i,j,k$ are distinct (4T)
  }\right. \right\rangle
\]
\noindent
Here, $4T$ is short for {\it four term relations}, and $L$ stands for 
{\it locality}.

One usually thinks of $\mathcal{A}_n$ in terms of horizontal chord diagrams on $n$
strands, where $t_{ij}$ is represented by a chord between strands $i$ and $j$,
and multiplication is done by stacking diagrams:

\vspace{2mm}
\begin{center}
\input{chorddiag.pstex_t}
\end{center}
\vspace{2mm}

Let us denote by $\mathcal{F}_2$ the free lie algebra over a field $k$ of 
characteristic 0, on two generators $X$ and 
$Y$, and by $\mathcal{UF}_2$ its universal enveloping algebra, which
is isomorphic to $k \langle \langle X,Y \rangle \rangle$, the algebra of non-commutative power series over $k$.
By a superscript $(m)$, as in $\mathcal{UF}_2^{(m)}$, we mean the space or object modulo
degree $(m+1)$, e.g. non-commutative polynomials of degree up to $m$.

There is a co-product on $\mathcal{UF}_2$, defined by $\Delta(X)=1 \otimes X+X \otimes 1$,
$\Delta(Y)=1 \otimes Y+Y \otimes 1$. An element $\Phi$ is called {\it group-like} if
$\Delta(\Phi)=\Phi \otimes \Phi$. An element $\varphi$ is
{\it primitive} (meaning it is a Lie-algebra element) if 
$\Delta(\varphi)=1 \otimes \varphi+\varphi \otimes 1$.

Now let us present the main characters of this story: the pentagon and hexagon
equations.

The pentagon equation originates from the fact that in a non-associative
algebra, there are five ways to multiply four elements, and these are connected by
a pentagon of re-associations. This fact and its parenthesized braid analogue are
shown in the figure below:
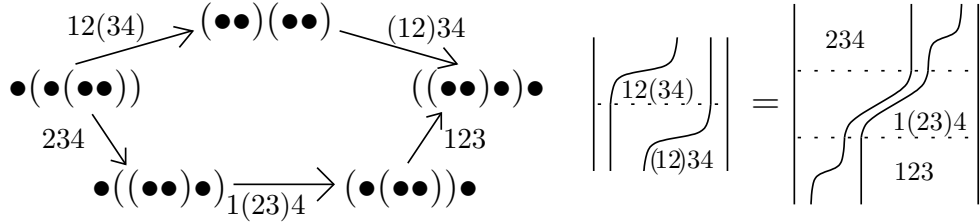
\begin{center}
 \input pentagon.pstex_t
\end{center}

The above equality of parenthesized braids implies an algebraic equation
in $\mathcal{A}_4$ called
the {\it pentagon equation}, $P(\Phi)=1$, for $\Phi \in {\mathcal UF}_2$, where
$$P(\Phi)=\Phi(t_{13}+t_{23},t_{34})^{-1}\Phi(t_{12},t_{23}+t_{24})^{-1}\Phi(t_{23},t_{34})
\Phi(t_{12}+t_{13},t_{24}+t_{34})\Phi(t_{12},t_{23}).$$

The hexagons arise from the following equivalences of parenthesized braids:

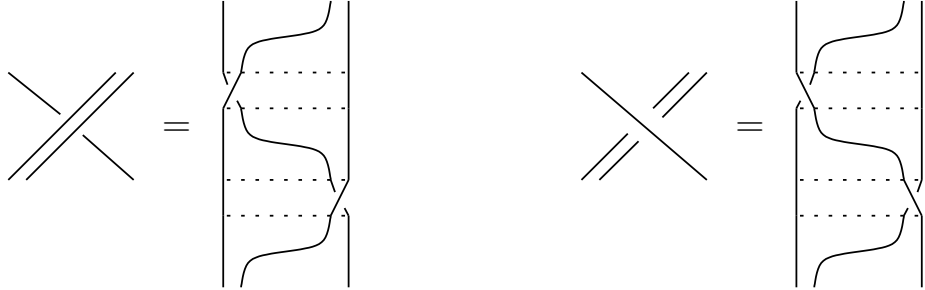
\begin{center}
 \input hexagon.pstex_t
\end{center}
Algebraically, the implied equations in $\mathcal{A}_3$ are $H_{\pm}(\Phi)=1$, for
$\Phi \in \mathcal{UF}_2$, where
$$H_{\pm}(\Phi)=e\left(\mp\frac{t_{13}+t_{23}}{2}\right)\Phi(t_{13},t_{12})e\left(\pm\frac{t_{13}}{2}\right)
\Phi(t_{13},t_{23})^{-1}e\left(\pm\frac{t_{23}}{2}\right)\Phi(t_{12},t_{23}).$$
Here,  $e(x):=e^x$.

An {\it associator} is a group-like element of $\mathcal{UF}_2$ which satisfies the 
pentagon and hexagon equations. 

Note that the spaces $\mathcal{A}_3$ and $\mathcal{UF}_2$ are almost isomorphic:
$\mathcal{UF}_2=\mathcal{A}_3/Z(\mathcal{A}_3)$ is the 
quotient of $\mathcal{A}_3$ by its center
$Z(\mathcal{A}_3)= \langle t_{12}+t_{23}+t_{31} \rangle$. An associator is sometimes defined
as an element of $\mathcal{A}_3$, which is equivalent, since associators are 
commutator-grouplike, meaning their abelianization is 1. (This follows
from the pentagon equation.)

\section{Theorem and Proof}

In \cite{F} Furusho proves the following surprising result:

\begin{theorem} \label{thm:main}
If $\Phi$ is a group-like element of $\mathcal{UF}_2$ with 
$c_2(\Phi)=\frac{1}{24}$, where $c_2(\Phi)$
denotes the coefficient of $XY$ in $\Phi$, 
and $\Phi$ satisfies the pentagon equation $P(\Phi)=1$ in $\mathcal{A}_4$,
then $\Phi$ satisfies the hexagon equations $H_{\pm}(\Phi)=1$
in ${\mathcal A}_3$, and therefore $\Phi$ is an associator. 
\end{theorem}

To prove Theorem~\ref{thm:main} we need two tools, which we 
will call {\it extension} and
{\it linearization} principles, and a major lemma, Lemma~\ref{lem:main} below. 

The ``extension principle''
is the assertion that an associator modulo degree $m$ can be lifted to
an associator modulo degree $(m+1)$. This statement is highly non-obvious; it takes up a
good part of \cite{Dr1} and the main part of \cite{BN}. Its main application is the
fact that non-trivial rational associators exist: indeed, 
up to degree $2$ a direct computation shows that $P(\Phi)=1$ and $H_{\pm}(\Phi)=1$
have a unique non-trivial solution (namely, $\Phi^{(2)}=1+\frac{1}{24}[X,Y]$), and
then by extension, one may construct a rational associator inductively.

The ``linearization principle'' is the standard fact that 
the pentagon $P(\Phi)=1$ and the hexagons
$H_{\pm}(\Phi)=1$ can be ``linearized''. Precisely, this means that there exist 
degree-preserving linear operators $dP: {\mathcal UF}_2 \to \mathcal{A}_4$ and 
$dH_{+}=dH_{-}=dH:\mathcal{UF}_2 \to \mathcal{A}_3$ so that if
$\Phi$ and $\Phi'$ satisfy the pentagon equation modulo degree $m$,
i.e. $P(\Phi)=P(\Phi')=1$ modulo degree $m$, and 
are equal modulo degree $m$, i.e. $\varphi:=(\Phi-\Phi')^{(m)}$ 
is homogeneous of degree $m$, then
\begin{center}
$P(\Phi)-P(\Phi') = dP(\varphi)$ modulo degree $(m+1)$ 
\end{center}
Likewise for the hexagons: if $\Phi$ and $\Phi'$ satisfy
one of the hexagons modulo degree $m$, and are equal modulo degree $m$,
then 
\begin{center}
$H_{\pm}(\Phi)-H_{\pm}(\Phi')=dH(\varphi)$ modulo degree $(m+1)$,
\end{center}
where $\varphi$ is as above. Note that if $\Phi$ and $\Phi'$ are group-like, then $\varphi$
is primitive.

To prove the theorem we use the following main lemma:
\begin{lemma} \label{lem:main} 
If $\varphi \in \mathcal{UF}_2$ is homogeneous of degree $m \geq 3$ and
is primitive, and satisfies
the linearized pentagon equation $dP(\varphi)=0$,
then $dH(\varphi)=0$, in other words the 
linearized pentagon equation implies the linearized hexagon equation.
\end{lemma}

\noindent{\bf Proof of Theorem \ref{thm:main} from Lemma \ref{lem:main}.} 
Let us assume that $\Phi$ is as in Theorem \ref{thm:main}, in particular
$P(\Phi)=1$, and, by contradiction, that one of the hexagons fails
to hold, i.e. $H_{\pm}(\Phi) \neq 1$. Let $m$ be the minimal degree in which 
this happens. 

By the simple computation in low degrees mentioned before, we know that 
$m \geq 3$. Note that by the minimality of $m$, $\Phi$
satisfies the hexagons modulo degree $m$ and hence it is 
an associator modulo degree $m$. By extension, there exists a 
$\Phi' \in \mathcal{UF}_2$ which agrees
with $\Phi$ modulo degree $m$, and which satisfies both the pentagon and
the hexagon equations modulo degree $(m+1)$. Let $\varphi=(\Phi-\Phi')^{(m)}$, 
which is homogeneous of degree $m$ since $\Phi=\Phi'$ modulo degree $m$. 

By linearization,
$dP(\varphi)=P(\Phi)-P(\Phi')=0$ modulo degree $(m+1)$, so by Lemma \ref{lem:main}, 
$H_{\pm}(\Phi)-H_{\pm}(\Phi')=dH(\varphi)=0$ modulo degree $(m+1)$. Back
again by linearization, as $\Phi'$ satisfies the hexagons in degree $m$,
it follows that so does $\Phi$, contradicting our pessimistic assumption from
the beginning. \qed

\section{Proof of the main lemma}
By explicit computations, 
\[dP(\varphi) =-\varphi(t_{12},t_{23}+t_{24})-\varphi(t_{13}+t_{23},t_{34})+\varphi(t_{23},t_{34})+
\varphi(t_{12}+t_{13},t_{24}+t_{34})+\varphi(t_{12},t_{23}),\]
\[dH(\varphi)=\varphi_(t_{13},t_{12})-\varphi(t_{13},t_{23})+\varphi(t_{12},t_{23}).\]

Let us start by proving two simple but necessary lemmas. Throughout this section, $\varphi$ is assumed to be a primitive, homogeneous
element of $\mathcal{UF}(2)$ of degree $\geq 3$, as in Lemma \ref{lem:main}. 

\begin{lemma}\label{as}
$dP(\varphi)=0$ implies that $\varphi$ is anti-symmetric, i.e.
$\varphi(X,Y)+\varphi(Y,X)=0.$
\end{lemma}

\parpic[r]{\input{qdef.pstex_t}}
\begin{proof}
We use the map $q: \mathcal{A}_4 \to \mathcal{F}_2$ defined by:
$t_{12} \mapsto X$, $t_{23} \mapsto Y$, $t_{13} \mapsto (-X-Y)$,
$t_{14} \mapsto Y$,  $t_{24} \mapsto (-X-Y)$, and $t_{34} \mapsto X$,
as illustrated by the figure on the right.
Since $q(dP(\varphi))=\varphi(X,Y)+\varphi(Y,X)$, the lemma follows.
\end{proof}

Note that by the anti-symmetry of $\varphi$, we have:
\[dP(\varphi)=\varphi(t_{12},t_{23})+\varphi(t_{34},t_{13}+t_{23})+
\varphi(t_{23}+t_{24},t_{12})+\varphi(t_{23},t_{34})+
\varphi(t_{12}+t_{13},t_{24}+t_{34}),\]
\[dH(\varphi)=\varphi(t_{12},t_{23})+\varphi(t_{23},t_{31})+\varphi(t_{31},t_{12}).\]

\begin{lemma}\label{eqlinhex}
The linearized hexagon
$dH(\varphi)=0$ is equivalent to the equation 
\begin{equation}\label{newdh}
\varphi(X,Y)+\varphi(Y,-X-Y)+\varphi(-X-Y,X)=0
\end{equation}
in $\mathcal{UF}_2$, if $\varphi$ is of degree $\geq 2$.
\end{lemma}

By an abuse of notation, we shall denote the left side of 
equation (\ref{newdh}) by $dH(X,Y)$.

\begin{proof}
$dH(\varphi)=0$ implies (\ref{newdh}) via the quotient map $\pi:\mathcal{A}_3 \to \mathcal{UF}_2$ 
which factors out by the center of $\mathcal{A}_3$, defined by $\pi(t_{12})=X$,
$\pi(t_{23})=Y$, and $\pi(t_{13})=-X-Y$.

For the other direction we apply 
the map $i: \mathcal{UF}_2 \to \mathcal{A}_3$
given by $i(X)=t_{12}$ and $i(Y)=t_{23}$

\noindent
$i(dH(X,Y))=\varphi(t_{12},t_{23})+\varphi(t_{23},-t_{12}-t_{23})+
\varphi(-t_{12}-t_{23},t_{12})=
\varphi(t_{12},t_{23})+
\varphi(t_{23},t_{13})+\varphi(t_{13},t_{12})=
dH(\varphi),$
where the second equality is due to the fact that $\varphi \in [\mathcal{F}_2, \mathcal{F}_2]$
(since it is primitive of degree $\geq 2$), and that the element 
$t_{12}+t_{23}+t_{13}$ is central in $\mathcal{A}_3$.
\end{proof}

{\bf Proof of the Main Lemma \ref{lem:main}}
Let us first introduce some notation.
Any map of sets $f: [r] \to [s]$, where $[r]:=\{1,2,...,r\}$ and $[s]=\{1,2,...,s\}$, induces
a map $\bar{f}: \mathcal{A}_s \to \mathcal{A}_r$, where
$\bar{f}(t_{ij}):=\sum_{\alpha \in f^{-1}(i), \beta \in f^{-1}(j)} t_{\alpha\beta}$.
If either of the pre-images is empty, we understand the sum to be zero.

Now, if $\psi(X,Y)\in\mathcal{UF}_2$, denote $\psi(123):=\psi(t_{12},t_{23}) \in 
\mathcal{A}_3$; and let 
\[\psi((i_1...i_t)(j_1...j_u)(k_1...k_v)):=\bar{g}(\psi(123)),\]
where $g:[t+u+v] \to [3]$ is given by $g^{-1}(1)=\{i_1...i_t\}$ and 
similarly for $2$ and $3$. 

The permutation group $S_4$ acts on $\mathcal{A}_4$ by commuting the
strands (i.e. the indices). For any $\sigma \in S_4$, we know that
$\sigma(dP(\varphi))=0$. We try to find a (small) set of permutations $\sigma_i$
such that $0=\sum _i \sigma_i dP(\varphi) = \sum_j \bar{g_j}(dH(\varphi))$,
for some $g_j:[4] \to [3]$. 

In fact, the following four permutations work 
(the notation $\sigma=x_1 x_2 x_3 x_4$ means 
$\sigma(i)=x_i$ for $i=1,2,3,4$):
\[\sigma_1=id, \quad \sigma_2=4231, \quad \sigma_3=1342, \quad \sigma_4=4312 \].
Using the notation introduced above,

\noindent
\begin{center}
$0=\sum_{i=1}^4 \sigma_i dP(\varphi)=$

$=\varphi(123)+
\cancel{\varphi(43(12))}+
\varphi((34)21)+\varphi(234)+\cancel{\varphi(1(23)4)}+$

$+\varphi(423)+
\cancel{\varphi(13(42))}+
\varphi((31)24)+\varphi(231)+\cancel{\varphi(4(23)1)}+$

$+\cancel{\varphi(134)}+
\varphi(24(13))+
\cancel{\varphi((42)31)}+
\varphi(342)+\varphi(1(34)2)+$

$+\cancel{\varphi(431)}+
\varphi(21(43))+
\cancel{\varphi((12)34)}+
\varphi(312)+\varphi(4(31)2)=$

$=dH(123)+dH((34)21)+dH(423)+dH((31)24),$
\end{center}
\vspace{2mm}

\noindent
where the cancellations are by the anti-symmetry of $\varphi$ (Lemma \ref{as}).

Note that every chord appearing on the right side ends on strand $2$,
i.e., $dH(123)+dH((34)21)+dH(423)+dH((31)24) \in \langle t_{12},t_{23},t_{24}\rangle .$
Also, $\langle t_{12},t_{23},t_{24}\rangle \cong \mathcal{F}_3 \subseteq \mathcal{A}_4,$
since there are no relations in $\mathcal{A}_4$ that involve only these elements.

Note that the anti-symmetry of $\varphi$ also implies the anti-symmetry of dH,
i.e. $dH(X,Y)=-dH(Y,X)$, and in particular $dH(X,X)=0$. We can now finish the
proof using two projections:

Let $p_1:\mathcal{F}_3 \to \mathcal{F}_2$ be the projection
defined by $t_{12} \mapsto X$, $t_{23} \mapsto Y$, $t_{24} \mapsto X$,
and apply this to the equality $0=dH(123)+dH((34)21)+dH(423)+dH((31)24)$.
We obtain
$0=dH(X,Y)+dH(X,Y)+dH(X+Y,X)+dH(X+Y,X),$
and therefore $dH(X+Y,X)=-dH(X,Y).$

Now do the same for the projection $p_2$ defined by $t_{12} \mapsto X$,
$t_{23} \mapsto X$, $t_{24} \mapsto Y$. We get
$0=dH(X,X)+dH(Y,X)+dH(X+Y,X)+dH(2X,Y)$, and so using the above,
we arrive to $dH(2X,Y)=2dH(X,Y)$.

This means that $dH(X,Y)$ contains only commutators that involve one $X$
and some number of $Y$'s, and so writing $dH$ in a linear basis
we have:

\begin{center}
$dH(X,Y)=\sum_{n=1}^\infty a_n(adY)^{n-1}(X).$
\end{center}

Since $dH(X,Y)=-dH(Y,X)$, we know that $a_n=0$ for all $n$
except possibly $n=2$, and $a_2=0$ by the assumption that
$\varphi$ is of degree $\geq 3$. Thus, $dH(X,Y)=0$, 
which is equivalent to the Main Lemma by Lemma \ref{eqlinhex}.
\qed

\vspace{3mm}
A surprising moment in the proof above is the equality relating a sum of four
linearized pentagon equations to a sum of linearized hexagons. As pointed out 
by Stavros Garoufalidis,  
there is a similar but more natural equality which, as the
reader can check, arises from the permuto-associahedron on page 21
of \cite{BN}:
$$dP(1234)-dP(1243)+dP(1423)-dP(4123)=dH(34(12))-dH((23)41)+dH(241)-dH(342)$$
It is easy to verify that the rest of the proof goes through the same way
using this equtlity in place of the one above.

\section{A side note}

The projection $q$, used in the proof of Lemma \ref{as},
has an interesting property and a nice topological interpretation.

If we embed $\mathcal{A}_3$ in $\mathcal{A}_4$ on any three strands 
(i.e. by any embedding of the index
set $[3]$ into $[4]$), and then apply $q$, we get the ``almost
isomorphism'' between $\mathcal{A}_3$ and $\mathcal{F}_2$, i.e. the
composition is factoring out by the center of $\mathcal{A}_3$:

\begin{center}
\input{qprop.pstex_t} 
\end{center}

This is a braid theoretic analogue of the following fact about the
symmetric group $S_4$. Since $S_4$ is isomorphic to the group of 
symmetries of the tetrahedron, and each element of $S_4$ also permutes
the three pairs of opposite edges of the tetrahedron, we obtain a map
$p: S_4 \to S_3$. Pre-composing this map with any embedding of $S_3$ into
$S_4$ induced by an embedding of the set $[3]$ into $[4]$,
we get an automorphism of $S_3$, namely, an isomorphism from
$S_3$ to the group of symmetries of a face:

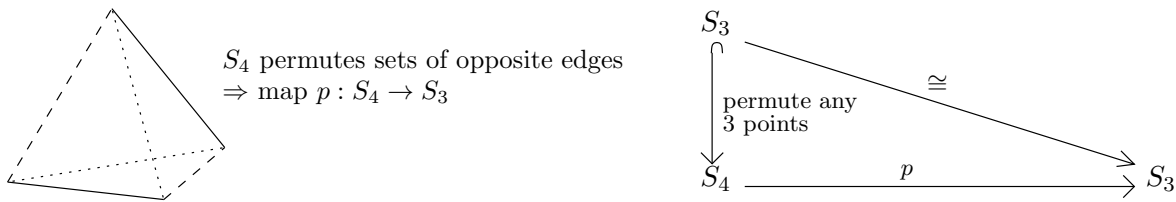
\begin{center}
\input pprop.pstex_t
\end{center}

Topologically, $q: \mathcal{A}_4 \to \mathcal{F}_2$ is induced by a
map $\tilde{q}$, defined as follows:

\begin{center}
 \input{topolq.pstex_t}
\end{center}

\parpic[r]{\begin{picture}(0,0)%
\includegraphics{spherical.pstex}%
\end{picture}%
\setlength{\unitlength}{3947sp}%
\begingroup\makeatletter\ifx\SetFigFont\undefined%
\gdef\SetFigFont#1#2#3#4#5{%
  \reset@font\fontsize{#1}{#2pt}%
  \fontfamily{#3}\fontseries{#4}\fontshape{#5}%
  \selectfont}%
\fi\endgroup%
\begin{picture}(988,558)(-3611,-7357)
\put(-3117,-7104){\makebox(0,0)[lb]{\smash{{\SetFigFont{7}{8.4}{\rmdefault}{\mddefault}{\updefault}{\color[rgb]{0,0,0}$=$}%
}}}}
\end{picture}%
}
Here, $pB_i$ denotes the pure braid group on $i$ strands, and $spB_4$
denotes the group of pure spherical braids on four strands. Spherical
braids live in $S^2 \times I$, as opposed to $D^2 \times I$, which means
that one strand wrapping all the way around the others is trivial, as
shown for strand $1$ on the right,
and similarly for all other strands. This defines the quotient map
$\tilde{q}_1$ above. 

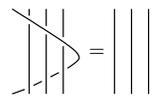
\parpic[r]{\begin{picture}(0,0)%
\includegraphics{twist.pstex}%
\end{picture}%
\setlength{\unitlength}{3947sp}%
\begingroup\makeatletter\ifx\SetFigFont\undefined%
\gdef\SetFigFont#1#2#3#4#5{%
  \reset@font\fontsize{#1}{#2pt}%
  \fontfamily{#3}\fontseries{#4}\fontshape{#5}%
  \selectfont}%
\fi\endgroup%
\begin{picture}(996,591)(-4308,-8608)
\put(-3837,-8340){\makebox(0,0)[lb]{\smash{{\SetFigFont{7}{8.4}{\rmdefault}{\mddefault}{\updefault}{\color[rgb]{0,0,0}$=$}%
}}}}
\end{picture}%
}
For $\tilde{q}_2$, take any spherical braid on
four strands, pull the last strand straight, and consider the first three
strands as a braid in the complement of strand $4$. The target space
of this map is the group of regular pure $3$-braids factored out by the
full twist of the three strands: pull the $4$-th strand straight on the
left side of the spherical relation shown in the figure on the right,
and observe that what we get is a full twist of the first $3$ strands,
which then has to be trivial in the image. 

Note that the chord diagram map induced by factoring out by the full twist in $pB_3$
is exactly the quotient map
$\pi: \mathcal{A}_3 \to \mathcal{F}_2$ which sends $t_{12}+t_{23}+t_{13}$ to $0$.

\end{document}